\documentclass{amsart}

\RequirePackage{color}

\def\smallskip{\vskip\smallskipamount}
\def\medskip{\vskip\medskipamount}
\def\bigskip{\vskip\bigskipamount}

\usepackage{amsmath,amsthm,bm,mathrsfs,amscd,tikz}
\usepackage{ytableau}
\usepackage{comment}
\usepackage{mathdots}
\usepackage[all]{xy}
\usepackage{graphicx}
\usepackage[colorinlistoftodos]{todonotes}
\usepackage{enumerate}
\usepackage{feynmp-auto}
\usepackage[english]{babel}
\usepackage[utf8]{inputenc}
\usepackage{float}
\usepackage{tikz}
\usepackage{tikz-cd}
\usepackage{amssymb}
\usepackage{mathtools}
\usetikzlibrary{decorations.pathmorphing}

\usepackage[utf8]{inputenc}

\newtheoremstyle{thmstyle}{}{}{\itshape}{}{\bfseries}{ }{5pt}{}
\newtheoremstyle{exstyle}{}{}{}{}{\bfseries}{ }{5pt}{}
\newtheoremstyle{defstyle}{}{}{}{}{\bfseries}{ }{5pt}{}
\newtheoremstyle{remstyle}{}{}{}{}{\bfseries}{ }{5pt}{}

\theoremstyle{thmstyle}
\newtheorem{thm}{Theorem}[section]
\newtheorem{theorem}[thm]{Theorem}

\newtheorem{lemma}[thm]{Lemma}

\newtheorem{proposition}[thm]{Proposition}

\newtheorem{corollary}[thm]{Corollary}

\newtheorem{conjecture}{Conjecture}[section]

\theoremstyle{exstyle}

\theoremstyle{defstyle}

\newtheorem{def-prop}[thm]{Definition-Proposition}
\newtheorem{def-lem}[thm]{Definition-Lemma}
\newtheorem{rem-convention}[thm]{Remark-Convention}
\newtheorem{def-note}[thm]{Definition-Notation}

\theoremstyle{remstyle}

\newtheorem{remark}[thm]{Remark}

\theoremstyle{remstyle}

\newcommand{\Hom}{\operatorname{Hom}}
\newcommand{\Fac}{\operatorname{Fac}}
\newcommand{\Ext}{\operatorname{Ext}}

\newcommand{\taurigid}{\operatorname{\tau-rigid}}

\DeclareMathOperator*{\rad}{rad}

\DeclareMathOperator*{\modu}{mod}
\DeclareMathOperator*{\Modu}{Mod}

\DeclareMathOperator*{\ann}{ann}
\DeclareMathOperator*{\Brick}{Brick}
\DeclareMathOperator*{\ind}{ind}
\DeclareMathOperator*{\Ind}{Ind}
\DeclareMathOperator*{\tors}{tors}
\DeclareMathOperator*{\Tors}{Tors}

\DeclareMathOperator*{\End}{End}

\DeclareMathOperator*{\brick}{brick}
\DeclareMathOperator*{\Ideal}{Ideal}
\DeclareMathOperator*{\rigid}{rigid}

\DeclareMathOperator*{\ttilt}{-tilt}
\DeclareMathOperator*{\tilt}{tilt}
\DeclareMathOperator*{\trig}{-rigid}

\DeclareMathOperator*{\soc}{soc}

\DeclareMathOperator*{\pd}{pd}

\DeclareMathOperator*{\Mri}{Mri}

\DeclareMathOperator*{\nD}{nD}

\DeclareMathOperator*{\sB}{sB}
\DeclareMathOperator*{\gC}{gC}

\DeclareMathOperator*{\GL}{GL}

\DeclareMathOperator*{\rep}{rep}

\makeatletter
\newcommand{\doublewidetilde}[1]{{%
  \mathpalette\double@widetilde{#1}%
}}
\newcommand{\double@widetilde}[2]{%
  \sbox\z@{$\m@th#1\widetilde{#2}$}%
  \ht\z@=.9\ht\z@
  \widetilde{\box\z@}%
}
\makeatother
\begin{document}

\title{Minimal ($\tau$-)tilting Infinite Algebras}
\author[Kaveh Mousavand, Charles Paquette]{Kaveh Mousavand, Charles Paquette} 
\address{Department of Mathematics and Statistics, Queen's University, Kingston ON, Canada}
\email{mousavand.kaveh@gmail.com}
\address{Department of Mathematics and Computer Science, Royal Military College of Canada, Kingston ON, Canada}
\email{charles.paquette.math@gmail.com}
\thanks{The first-named author was partially supported by Mitcas Globalink Research Award. The second-named author was supported by the National Sciences and Engineering Research Council of Canada, and by the Canadian Defence Academy Research Programme.}


\maketitle

\begin{abstract}
Motivated by a new conjecture on the behavior of bricks, we start a systematic study of minimal $\tau$-tilting infinite algebras. In particular, we treat minimal $\tau$-tilting infinite algebras as a modern counterpart of minimal representation infinite algebras and show some of the fundamental similarities and differences between these families. We then relate our studies to the classical tilting theory and observe that this modern approach can provide fresh impetus to the study of some old problems. We further show that in order to verify the conjecture it is sufficient to treat those minimal $\tau$-tilting infinite algebras where almost all bricks are faithful.
Finally, we also prove that minimal extending bricks have open orbits, and consequently obtain a simple proof of the brick analogue of the First Brauer-Thrall Conjecture, recently shown by Schroll and Treffinger using some different techniques.
\end{abstract}

\tableofcontents

\section{Introduction}\label{Section:Introduction}

Throughout this note, $k$ denotes an algebraically closed field and all algebras are assumed to be finite dimensional over $k$, associative and unital. 
We denote by $\Modu \Lambda$ the category of all left $\Lambda$-modules, and $\modu \Lambda$ denotes the full subcategory of finitely generated left $\Lambda$-modules.
Let $\Ind(\Lambda)$ and $\ind(\Lambda)$ respectively denote the collections of all indecomposable modules in $\Modu \Lambda$ and $\modu \Lambda$, up to isomorphism. Then, $\Lambda$ is said to be \emph{representation infinite} if $\ind(\Lambda)$ is infinite. In this work we are primarily interested in representation infinite algebras, because the problems treated here are trivial over representation finite algebras.

For our purposes in this paper, without loss of generality, we may always assume that $\Lambda$ is basic and connected. Thus, there is an isomorphism of algebras $\Lambda \simeq kQ/I$, where $Q$ is a finite connected quiver and $I$ is an admissible ideal in the path algebra $kQ$. In particular, $|K_0(\Lambda)|=|Q_0|$, meaning that the rank of the Grothendieck group $K_0(\Lambda)$ is the same as the number of vertices in $Q$.
Consequently, each $\Lambda$-module can be viewed as a representation of the bound quiver $(Q,I)$.
For the rudiments of representation theory of algebras we refer to \cite{ASS}.

\subsection{Motivations and background}
In 2014, Adachi, Iyama and Reiten \cite{AIR} introduced $\tau$-tilting theory of associative algebras as a modern generalization of classical tilting theory.
For an algebra $\Lambda$, $\tau$-tilting theory treats the notion of rigidity and compatibility of $\Lambda$-modules with respect to the Auslander-Reiten translation $\tau_{\Lambda}$, which we often denote by $\tau$.
In particular, the set of (support) $\tau$-tilting modules extends that of the classical tilting modules and admits a more consistent behavior with respect to the notion of mutation (see \cite{AIR} for details).

For an algebra $\Lambda$, let $\tilt(\Lambda)$ denotes the set of all isomorphism classes of basic tilting modules in $\modu \Lambda$. Similarly, by $\tau \ttilt(\Lambda)$ we denote the set of all basic $\tau$-tilting modules in $\modu \Lambda$, up to isomorphism. 
Then, $\Lambda$ is called \emph{tilting finite} if $|\tilt(\Lambda)| < \infty$ and it is \emph{$\tau$-tilting finite} if $|\tau \ttilt(\Lambda)| < \infty$.
As shown in \cite{AIR}, (support) $\tau$-tilting modules closely relate to several fundamental objects in $\modu \Lambda$ and its bounded derived category.
Thus, for a given algebra, it is important to decide whether it is $\tau$-tilting (in)finite.
Since $\tilt(\Lambda) \subseteq \tau\ttilt(\Lambda)$, each $\tau$-tilting finite algebra is evidently tilting finite, but the converse does not hold in general.

Recall that a module $M$ in $\Ind(\Lambda)$ is a \emph{brick} if $\End_{\Lambda}(M)$ is a division algebra. Suppose $\Brick(\Lambda)$ denotes a collection of all  bricks in $\Modu \Lambda$, up to isomorphism, and similarly, $\brick(\Lambda)$ is a set of all bricks in $\modu \Lambda$, up to isomorphism.
In \cite{DIJ}, the authors established the ``$\tau$-rigid-brick correspondence", which is an elegant linkage between those indecomposable modules $M$ in $\modu \Lambda$ for which $\Hom_{\Lambda}(M,\tau M)=0$, and elements of $\brick(\Lambda)$.
Then, they concluded that $\Lambda$ is $\tau$-tilting finite if and only if $\brick(\Lambda)$ is finite. More importantly, from \cite{DIJ} we know that $\tau$-tilting finite algebras can be viewed as a natural generalization of representation-finite algebras and they admit a phenomenon similar to the finite type cluster algebras: for $\Lambda=kQ/I$ with $|Q_0|=n$, each vertex of the mutation graph of support $\tau$-tilting modules in $\modu \Lambda$ is $n$-regular (i.e, each vertex is of degree $n$) and the graph is finite if and only if $\Lambda$ is $\tau$-tilting finite.

Notice that the problem of $\tau$-tilting finiteness is of interest only at the level of representation infinite (rep-inf, for short) algebras. Therefore, it is natural to treat $\tau$-tilting finiteness of those algebras which are rep-inf and minimal with respect to this property. Recall that an algebra $\Lambda$ is said to be \emph{minimal representation-infinite} (min-rep-inf, for short) if $\ind(\Lambda)$ is infinite but $\ind(\Lambda/J)$ is a finite set, for every non-zero ideal $J$ in $\Lambda$. 
For any rep-inf algebra $\Lambda$, let $\Mri(\Lambda)$ be the set of isomorphism classes of all quotient algebras of $\Lambda$ that are min-rep-inf.
Observe that if there exists $\Lambda'$ in $\Mri(\Lambda)$ which is $\tau$-tilting infinite, then so is $\Lambda$. Thus, it is important to have a description of $\tau$-tilting (in)finiteness of min-rep-inf algebras. This direction of work has been developed by the first author in \cite{Mo1, Mo2}. For further details on this classsification, see subsection \ref{Subsection:Minimality conditions}.

Although $\tau$-tilting finite algebras have been treated extensively, the study of $\tau$-tilting infinite algebras is significantly more complicated and among them one can hope for a better description of those which are minimal with respect to this property. We say $\Lambda$ is \emph{minimal $\tau$-tilting infinite} (or min-$\tau$-infinite, for short) if $\Lambda$ is a $\tau$-tilting infinite algebra but $\Lambda/J$ is $\tau$-tilting finite, for each nonzero idea $J$ in $\Lambda$. 
Our interest in min-$\tau$-infinite algebras primarily stems from their decisive role in the study of the following modern conjecture.

\begin{conjecture}\label{main conjecture}
If $\Lambda$ is $\tau$-tilting infinite, there exists a positive integer $d$ such that there are infinitely many (non-isomorphic) bricks of length $d$ in $\modu \Lambda$.
\end{conjecture}

The above conjecture can be viewed as a modern analogue of the celebrated Brauer-Thrall conjectures, yet it is not a verbatim counterpart of any of them in terms of bricks. 
In particular, the number of integers $d$ in the assertion of Conjecture \ref{main conjecture} can be finite or infinite (in the $2$-Kronecker algebra, there is exactly one such $d$, whereas for the the $3$-Kronecker case we get infinitely many such $d$). 

Conjecture \ref{main conjecture} first appeared in \cite{Mo2}, where the first-named author phrased it in a more geometric language to establish a new linkage between $\tau$-tilting theory and geometry of representation varieties of quivers.
Further, in the aforementioned paper the conjecture was verified for gentle algebras, as well as all the min-rep-inf algebras treated there. In \cite{STV}, the authors have also come across the same conjecture and verified it for the family of special biserial algebras.

To prove Conjecture \ref{main conjecture} in full generality, we observe that it suffices to show it for all min-$\tau$-infinite algebras. Further, note that via the notion of generic modules and their properties, we can also view the above conjecture from a new perspective. In particular, in \cite{C-B} Crawley-Boevey shows that $\Lambda$ is representation infinite if and only if for some $d$, there are infinitely many (non-isomorphic) modules in $\ind(\Lambda)$ of endolength $d$. 
Adopting a similar viewpoint, Conjecture \ref{main conjecture} asserts that ``$\Lambda$ is brick-infinite if and only if for some $d$, there are infinitely many (non-isomorphic) modules in $\brick(A)$ of (endo)length $d$".

\subsection{Outline and main results}
In Section \ref{Section:Preliminaries}, we introduce our notations and terminology. Further, we collect some known results essential in our work. 
Before we state our main results, let us give analogous characterizations of the classical and modern minimality conditions discussed in Section \ref{Section:Preliminaries}.
The following statements follow from a known result of Auslander \cite{Au} and the recent work of Sentieri \cite{Se}. In particular, for an algebra $\Lambda$, we have

\medskip

\begin{itemize}
    \item $\Lambda$ is minimal representation infinite  if and only if $\Ind(\Lambda)\setminus \ind(\Lambda)\neq \emptyset$ and every $M \in \Ind(\Lambda)\setminus \ind(\Lambda)$ is faithful.
    
    \item $\Lambda$ is minimal $\tau$-tilting  infinite if and only if $\Brick(\Lambda)\setminus \brick(\Lambda)\neq \emptyset$ and every $N \in \Brick(\Lambda)\setminus \brick(\Lambda)$ is faithful.
\end{itemize}

\medskip

The next theorem lists some important properties of min-$\tau$-infinite algebras. As explained in Section \ref{Section: Minimal tau-tilting infinite algebras}, this theorem also highlights some fundamental differences and similarities between the min-$\tau$-infinite algebras and the min-rep-inf algebras.

\begin{theorem}
Let $\Lambda=kQ/I$ be a minimal $\tau$-tilting infinite algebra. Then $\Lambda$ is node-free and central. Moreover, $\Lambda$ admits no projective-injective module.
\end{theorem}

The previous theorem leads to the following classification result on algebras with radical square zero, from which one immediately verifies Conjecture \ref{main conjecture} for such algebras.
This also gives a simple proof for the following result of Adachi \cite{Ad}.
We say $Q$ is a \emph{sink-source} quiver if every vertex of $Q$ is either a sink or a source.

\begin{corollary}[\cite{Ad}]
Let $\Lambda=kQ/I$ be such that $\rad^2(\Lambda)=0$. The following are equivalent:
\begin{enumerate}
    \item $\Lambda$ is $\tau$-tilting infinite;
    \item $Q$ contains a sink-source subquiver of affine type.
\end{enumerate}
\end{corollary}

In Section \ref{Section:Minimal tilting infinite algebras} we slightly change our perspective and study the minimal $\tau$-tilting infinite algebras from the viewpoint of classical tilting theory. 
More specifically, we prove that minimal $\tau$-tilting infinite algebras are in fact minimal with respect to the notion of tilting modules, where the minimality is defined analogously: An algebra $\Lambda$ is \emph{minimal tilting infinite} if it admits infinitely many tilting modules, but every proper quotient algebra $\Lambda/J$ has only finitely many tilting modules.

Before we summarize the main results of Section \ref{Section:Minimal tilting infinite algebras} in the following theorem, let us recall that for a collection of objects $\mathfrak{O}$, we say \emph{almost all} objects of $\mathfrak{O}$ satisfy property $\mathcal{P}$ provided all but finitely many objects of $\mathfrak{O}$ have property $\mathcal{P}$.

\begin{theorem}
Let $\Lambda$ be an algebra. If $\Lambda$ is minimal $\tau$-tilting infinite, then almost all $\tau$-rigid $\Lambda$-modules are partial tilting. Consequently, we have:
\begin{enumerate}
    \item $\Lambda$ is minimal $\tau$-tilting infinite if and only if it is minimal tilting infinite.
    \item If $\Lambda$ is minimal $\tau$-tilting infinite, then the mutation graph of tilting modules in $\modu \Lambda$ is infinite and regular at almost every vertex.
\end{enumerate}
\end{theorem}

In Section \ref{Section:Reduction of the conjecture} we use minimal $\tau$-tilting infinite algebras in the study of Conjecture \ref{main conjecture}.
Before we state the next theorem, we note that a min-$\tau$-infinite algebra may admit infinitely many non-isomorphic bricks with distinct annihilators. 
However, we prove that it is sufficient to only treat the conjecture for a particular subfamily of min-$\tau$-infinite algebras, as stated in the following theorem.

\begin{theorem}\label{Thm: reduction theorem in introduction}
To verify Conjecture \ref{main conjecture}, it is sufficient to consider minimal $\tau$-tilting infinite algebras for which almost all bricks are faithful.
\end{theorem}

Finally, in Section \ref{Section: First BT conjecture} we employ some tools from the algebro-geometric aspects of representation theory of algebras to show the following result on the minimal extending bricks (for the definition, see Section \ref{Section:Reduction of the conjecture}). As an immediate consequence of this theorem, we give a simple proof of the modern analogue of the First Brauer-Thrall Conjecture, first shown in \cite{ST}.

\begin{theorem}
Let $\Lambda$ be an algebra and $M$ a minimal extending brick in $\modu \Lambda$ with dimensional vector $\bf{d}$. Then, the orbit of $M$ under the $\GL(\bf{d})$ action is open.
\end{theorem}

As mentioned before, our work is primarily motivated by Conjecture \ref{main conjecture} and some of our main arguments rely on the study of bricks. Nevertheless, we stated our results in terms of $\tau$-tilting (in)finineness of algebras under consideration. This allows us to better relate our study to the other subjects discussed in the paper.

\section{Preliminaries}\label{Section:Preliminaries}
In this section we only collect some basic materials needed in the rest of the paper.
The well-know results appear without proofs but only references, if necessary.

\subsection{Notations and conventions}
In addition to those introduced in Section \ref{Section:Introduction}, here we fix some notations and conventions used throughout the paper.

By a quiver we always mean a finite directed graph, formally given by a quadruple $Q=(Q_0,Q_1,s,e)$, with the vertex set $Q_0$ and arrow set $Q_1$, and the functions $s,e: Q_1 \rightarrow Q_0$ respectively send each arrow $\alpha$ to its start $s(\alpha)$ and its end $e(\alpha)$.
As mentioned before, up to Morita equivalence, every algebra $\Lambda$ treated in this paper is of the form $\Lambda=kQ/I$, for a unique quiver $Q$ and an admissible ideal $I$ in $kQ$.
In particular, $|K_0(\Lambda)|=|Q_0|$, where $|Q_0|$ is the number of vertices of $Q$.
In a bound quiver $(Q,I)$, if the arrow $\alpha$ ends where the arrow $\beta$ starts, then $\beta \alpha$ denotes the path of length two in $Q$. 
Moreover, a vertex $v$ in $(Q,I)$ is called a \emph{node} if it is neither a sink nor a source, and for each pair of arrows $\alpha$ and $\beta$ in $Q$ with $e(\alpha)=v=s(\beta)$, we have $\beta\alpha \in I$. Then, $\Lambda=kQ/I$ is called \emph{node-free} if $(Q,I)$ has no nodes. 

If $\Lambda=kQ/I$, by $R$ we denote a minimal set of (uniform) relations that generates $I$. Namely, each element of $R$ is a linear combination of the form $r=\sum_{i=1}^{t} \lambda_i p_i$, with $t\in \mathbb{Z}_{>0}$ and $\lambda_i \in k\setminus \{0\}$, and every $p_i$ is a path in $Q$ whose length is not smaller than two, and all $p_i$ start at the same vertex and also end at the same vertex. The relation $r$ is \emph{monomial} if $t=1$, and is called \emph{binomial} if $t=2$.

For any algebra $\Lambda$, let $\Ideal(\Lambda)$ denote the set of all (two-sided) ideals in $\Lambda$. Note that $\Ideal(\Lambda)$ carries a natural lattice structure, where $J_1 \wedge J_2:=J_1 \cap J_2$ and $J_1 \vee J_2:=J_1 + J_2$, for each pair $J_1$ and $J_2$ in $\Ideal(\Lambda)$. Consequently, $\Lambda$ is called \emph{distributive} if $\Ideal(\Lambda)$ is a distributive lattice. Otherwise, $\Lambda$ is \emph{non-distributive}. It is immediate that every quotient of a distributive algebra is again distributive. Moreover, it is known that each non-distributive algebra is representation infinite. 
We also recall that a $k$-algebra is said to be \emph{central} if $Z(\Lambda)=k$, where $Z(\Lambda)$ denotes the center of $\Lambda$.

For $X$ in $\modu \Lambda$, let $|X|$ denote the number of non-isomorphic indecomposable summands of $X$. Then, $X$ is \emph{basic} if $|X|$ is exactly the same as the number of indecomposable summands of $X$.
Moreover, by $\pd_{\Lambda}(X)$ and $\tau_{\Lambda}X$ we respectively denote the projective dimension and the Auslander-Reiten translation of $X$ in $\modu \Lambda$. If there is no confusion, $\Lambda$ is often suppressed from our notations.
Support of $X$ at vertex $i\in Q_0$ is defined as  $X_i:=\Hom_{\Lambda}(P_i,X)$, and the dimension vector of $X$ is  $\underline{\dim}(X)=(\dim X_i)_{i\in Q_0}$.
Then, $X$ is \emph{sincere} if $\underline{\dim}(X) \in \mathbb{Z}^{Q_0}_{>0}$.

An algebra $\Lambda$ is \emph{biserial} if for every left and right indecomposable projective $\Lambda$-module $P$, the submodule $\rad(P)$ is sum of at
most two uniserial modules $X$ and $Y$ and $X \cap Y$ is either zero or a simple module. Moreover, $\Lambda$ is \emph{special biserial} if it is Morita equivalent to an algebra $kQ/I$ such that $(Q,I)$ satisfies the following:
\begin{enumerate}[(B1)]
\item At every vertex $x$ in $Q$,  there are at most two incoming and at most two outgoing arrows.
\item For each arrow $\alpha$, there is at most one arrow $\beta$ such that $\beta \alpha \notin I$ and at most one arrow $\gamma$ such that $\alpha \gamma \notin I$.
\end{enumerate}

We note that every special biserial algebras is biserial, but the converse is not true. Further, for each special biserial algebra $\Lambda=kQ/I$, one can choose a minimal set of generators for $I$ consisting of only monomial and binomial relations.
A special biserial algebra $\Lambda=kQ/I$ is said to be a \emph{string} algebra if there is a set of monomial relations in $kQ$ which generate $I$.

By a subcategory of $\Modu \Lambda$, we always mean a full subcategory which is closed under direct sum, direct summands, and isomorphisms.
A subcategory $\mathcal{T}$ of $\modu \Lambda$ is a \emph{torsion class} if it is closed under quotients and extensions. Dually, a subcategory $\mathcal{F}$ of $\modu \Lambda$ is \emph{torsion free} if $\mathcal{F}$ is closed under submodules and extensions. It is well-known that every torsion class $\mathcal{T}$ uniquely determines a torsion pair $(\mathcal{T}, \mathcal{F})$ in $\modu \Lambda$, where for each $X \in \mathcal{T}$ and every $Y \in \mathcal{F}$, we have $\Hom_{\Lambda}(X,Y)=0$, and $\mathcal{F}$ is maximal with this property. 

\subsection{$\tau$-tilting finite algebras}

A $\Lambda$-module $M$ is called \emph{rigid} if $\operatorname{Ext}^1_{\Lambda}(M,M)=0$. By $\rigid (\Lambda)$ we denote the set of all basic rigid objects in $\modu \Lambda$, considered up to isomorphisms. 
Similarly, $M$ is said to be \emph{$\tau$-rigid} if $\Hom_{\Lambda}(M,\tau M)=0$. Analogously, we use $\tau \trig(\Lambda)$ and $\mathtt{i}\tau \trig(\Lambda)$ to respectively denote the set of basic $\tau$-rigid modules and the indecomposable $\tau$-rigid objects in $\modu \Lambda$, up to isomorphisms.
A rigid module $X$ is called \emph{tilting} if $\pd_{\Lambda}\leq 1$ and $|X|=|\Lambda|$. 
Similarly, each $\tau$-rigid module $M$ with $|M|=|\Lambda|$ is said to be \emph{$\tau$-tilting}. More generally, $M$ is called \emph{support $\tau$-tilting} if $M$ is $\tau$-tilting over $\Lambda/\langle e \rangle$, where $e$ is a maximal idempotent in $\ann_{\Lambda}(M)$. 
By $\tau \ttilt(\Lambda)$ and $s\tau \ttilt(\Lambda)$ we respectively denote the set of basic $\tau$-tilting modules and that of basic support $\tau$-tilting modules in $\modu \Lambda$, considered up to isomorphisms.

Recently there have been various attempts to systematically study the $\tau$-tilting finiteness of algebras. That is to find the necessary and sufficient conditions for arbitrary algebra $\Lambda$ such that $|\tau \ttilt(\Lambda)| < \infty$. 
The ``brick-$\tau$-rigid correspondence" introduced by Demonet, Iyama and Jasso \cite{DIJ} was a significant step towards such an objective. 
The following theorem collects some fundamental characterizations of $\tau$-tilting finite algebras extensively used in this paper.

\begin{theorem}[\cite{AIR, DIJ}]\label{Thm:tau-tilting finiteness equivalences}
For an algebra $\Lambda$, the following are equivalent:
\begin{enumerate}
\item  $\Lambda$ is $\tau$-tilting finite;
\item  $s\tau \ttilt(\Lambda)$ is finite;
\item $\tau \rigid(\Lambda)$ is finite;
\item $\modu \Lambda$ contains only finitely many isomorphism classes of bricks;
\item  Every torsion(-free) class in $\modu \Lambda$ is functorially finite;
\end{enumerate}
\end{theorem}

\subsection{Minimality conditions}\label{Subsection:Minimality conditions}
Before we review some minimality conditions, let us make a handy observation that will be freely used throughout our work.

We recall that any epimorphism of algebras $\psi: \Lambda_1 \rightarrow \Lambda_2$ induces an exact functorial embedding $\widetilde{\psi}: \modu\Lambda_2 \rightarrow \modu\Lambda_1$. Particularly, every indecomposable (respectively, brick) in $\modu \Lambda_2$ can be seen as an indecomposable (respectively, brick) in $\modu \Lambda_1$.
Hence, representation finiteness of algebras is preserved under algebraic quotients: if $\Lambda$ is rep-finite, then $\Lambda/J$ is rep-finite, for each $J$ in $\Ideal(\Lambda)$.
Moreover, by Theorem \ref{Thm:tau-tilting finiteness equivalences}, it is immediate that the $\tau$-tilting finiteness of algebras is also preserved under taking quotients. 
In constrast, we note that there are tilting finite algebras $\Lambda$ such that $\Lambda/J$ is tilting infinite, for some $J \in \Ideal (\Lambda)$.

Due to the above observation, representation infinite algebras which are minimal with respect to this property have been decisive in the study of representation theory of algebras, such as in the proofs of Brauer-Thrall conjectures.
As defined earlier, $\Lambda$ is \emph{minimal representation infinite} (or min-rep-inf, for short) if it is representation infinite but any proper quotient of $\Lambda$ is representation finite.
The literature on min-rep-inf algebras is very rich (for instance, see \cite{Bo1, Bo2} and the references therein). A conceptual classification of these algebras appears in \cite{Ri}, where Ringel shows that each min-rep-inf algebra belongs to at least one of the following families. For the definition of a \emph{good covering}, see \cite{Bo1}.

\begin{itemize}
\item $\Mri({\mathfrak{F}_{\sB}})$: min-rep-inf \textbf{s}pecial \textbf{B}iserial algebras;

\item $\Mri({\mathfrak{F}_{\nD}})$: min-rep-inf \textbf{n}on-\textbf{D}istributive algebras;

\item $\Mri({\mathfrak{F}_{\gC}})$: min-rep-inf algebras with a \textbf{g}ood \textbf{C}overing $\widetilde{\Lambda}$ such that a finite convex subcategory of $\widetilde{\Lambda}$ is tame-concealed of type $\widetilde{\mathbb{D}}_n$ or $\widetilde{\mathbb{E}}_{6,7,8}$.
\end{itemize}

The bound quivers of the first two subfamilies from the above list are respectively described in \cite{Ri} and \cite{Bo2}. Moreover, it is known that a concrete classification of the last subfamily in terms of their bound quivers is very hard, if not impossible, and soon falls out of control as the number of simple modules grows. 

As of a more contemporary approach to the study of min-rep-inf algebras, it is natural to determine which ones are $\tau$-tilting (in)finite.
One can show that each of the above-mentioned subfamilies contains both $\tau$-tilting finite, as well as $\tau$-tilting infinite algebras. 
In \cite{Mo1} and \cite{Mo2}, the first-named author carried out a full study of $\tau$-tilting finiteness of the algebras in  $\Mri({\mathfrak{F}_{\sB}})$ and $\Mri({\mathfrak{F}_{\nD}})$ and determined which ones are $\tau$-tilting finite and which ones are not. 
The following theorem captures the main results in this direction. For explicit description of the bound quivers and further details, see the aforementioned papers.

\begin{theorem}[\cite{Mo2}]
Let $\Lambda= kQ/I$ belong to $\Mri({\mathfrak{F}_{\sB}}) \cup \Mri({\mathfrak{F}_{\nD}})$. Then $\Lambda$ is $\tau$-tilting finite if and only if $(Q,I)$ has a node or a non-quadratic monomial relation.
\end{theorem}

\section{Minimal $\tau$-tilting infinite algebras}\label{Section: Minimal tau-tilting infinite algebras}

As already noticed, $\tau$-tilting finiteness of algebras is preserved under taking algebraic quotients. 
In \cite{Mo1} and \cite{Wa}, where the problem of $\tau$-tilting finiteness is concerned, the authors independently noticed that a natural analogue of minimal representation infinite algebras in the modern setting should be crucial in the study of $\tau$-tilting (in)finiteness. In particular, an algebra $\Lambda$ is \emph{minimal $\tau$-tilting infinite} (or min-$\tau$-infinite, for short), if $\Lambda$ is $\tau$-tilting infinite but every proper quotient algebra of $\Lambda$ is $\tau$-tilting finite.

Inspired by the classification of min-rep-inf algebras, as briefly recalled in Section \ref{Subsection:Minimality conditions}, it is natural to search for a classification of min-$\tau$-infinite algebras. Although it might be impossible to describe the bound quivers of all of min-$\tau$-infinite algebras, one can still treat this problem for certain subfamilies of them. We pursue this direction of work in a separate paper.
In this section we focus on some fundamental properties of min-$\tau$-infinite algebras. 
In particular, we compare min-rep-inf and min-$\tau$-infinite algebras and derive some interesting results from this comparison.

It is known that a min-rep-inf algebra admits no projective-injective module. Let us begin by showing that the same property holds for min-$\tau$-infinite algebras. For the sake of brevity, for two algebras $\Lambda_1$ and $\Lambda_2$, we say $\Lambda_1$ and $\Lambda_2$ are of the same \emph{$\tau$-tilting type} when $\Lambda_1$ is $\tau$-tilting finite if and only if $\Lambda_2$ is so.

\begin{proposition}
Let $\Lambda$ be an algebra such that $\modu \Lambda$ contains a projective-injective module. Then, there exists a nonzero ideal $J$ such that $\Lambda$ and $\Lambda/J$ are of the same $\tau$-tilting type. 
Consequently, if $\Lambda$ is minimal $\tau$-tilting infinite, there is no projective-injective module in $\modu \Lambda$.
\end{proposition}
\begin{proof}
Suppose $\Lambda=kQ/I$ and let $X \in \ind(\Lambda)$ be projective-injective. Then, $X=P_x=I_z$, for a pair of vertices $x$ and $z$ in $Q_0$.
Hence, $X= \Lambda e_x$ and $\soc(X)$ is a subspace of $e_z \Lambda e_x$. 
Note that $\soc(X)$ is simple and $\Lambda$ is basic. Then we have $\dim_k(\soc(X))=1$.
Consider the ideal $J:=\langle \rho \rangle$, where $\rho$ denotes the element of $e_z \Lambda e_x \subseteq kQ/I$ associated to $\soc(X)$. Then $J$ is one-dimensional and $JY=0$ if and only if $Y\not \simeq X$.
Consequently, $\brick(\Lambda) \setminus \{X\} \subseteq \brick (\Lambda/J) \subseteq \brick(\Lambda)$. Then, the first part of the assertion follows from Theorem \ref{Thm:tau-tilting finiteness equivalences}.
Moreover, the second part is an immediate consequence of the minimality assumption.
\end{proof}

For $\Lambda=kQ/I$, suppose $a$ and $z$ are respectively a source and a sink in $(Q,I)$. Let $(Q',I')$ be the bound quiver obtained from $(Q,I)$ by \emph{gluing} the vertices $a$ and $z$, as follows: First identify $a$ and $z$, then kill all the composition of arrows $\beta \alpha$, for each $\alpha$ incoming to $z$ with any $\beta$ outgoing from $a$. 
If $v$ denotes a vertex of $(Q',I')$ obtained from gluing a sink and a source in $(Q,I)$, then obviously $v$ is a node. The reverse process of gluing is called \emph{resolving} a node.
It is well-known that if $\Lambda'$ is obtained from $\Lambda$ via a sequence of gluing (resolving), then $\Lambda$ and $\Lambda'$ are of the same representation type. In particular, $\Lambda$ is min-rep-inf algebra if and only if $\Lambda'$ is so.

The following proposition shows that there is a fundamental difference between the bound quivers of min-$\tau$-infinite algebras and min-rep-inf algebras. This result is already shown in \cite{Mo2}, thus we omit the proof.

\begin{proposition}[\cite{Mo2}]\label{Prop: No node property}
Let $\Lambda=kQ/I$ be a minimal $\tau$-tilting infinite algebra. Then, $(Q,I)$ does not contain any node.
\end{proposition}

By the above proposition and the remark preceding that, one observes that for the classification of min-rep-inf algebras and min-$\tau$-infinite algebras we only need to treat the node-free bound quivers. Moreover, the proposition gives a simple sufficient condition for $\tau$-tilting finiteness of a large family of min-rep-inf algebras. We note that, however, there also exist node-free min-rep-inf algebras which are $\tau$-tilting finite (for explicit examples, see the wind wheel algebras in \cite{Mo1}). 

As an important consequence of the previous proposition, we get an explicit criterion for $\tau$-tilting (in)finiteness of algebras with radical square zero. 
We note that a classification of $\tau$-tilting finite algebras with radical square zero also appears in \cite{Ad}, where Adachi uses the notion of separated quivers.
However, as shown in our proof, the reduction to the min-$\tau$-infinite algebras allows us to use the property of their bound quivers and obtain the same result via a less technical approach.

To state the next result more succinctly, we say $Q$ is a \emph{sink-source quiver} if each vertex of $Q$ is either a sink or a source.
\begin{corollary}[\cite{Ad}]\label{Thm: radical square zero}
Provided that $\rad^2(\Lambda)=0$, the following are equivalent:
\begin{enumerate}
    \item $\Lambda$ is $\tau$-tilting infinite;
    \item There exists an ideal $J$ in $\Lambda$ such that $\Lambda/J$ is hereditary of affine type with a sink-source quiver.  
\end{enumerate}
\end{corollary}
\begin{proof}
Suppose $\Lambda$ is $\tau$-tilting infinite. Without loss of generality, we can assume $\Lambda$ is min-$\tau$-infinite. 
By Proposition \ref{Prop: No node property}, $\Lambda$ is node-free. Therefore, $\rad^2(\Lambda)=0$ implies that $\Lambda=kQ$, for a sink-source quiver $Q$. To finish the proof, we use the fact that hereditary algebras of affine type are $\tau$-tilting infinite and obviously minimal with respect to this property. 
\end{proof}

For an affine quiver $Q$, it is well-known that the path algebra $kQ$ admits a one-parameter family of bricks of the same length (for example, see \cite{SS}). Therefore, the following result is an immediate consequence of the above classification and verifies Conjecture \ref{main conjecture} for the family of algebras with $\rad^2(\Lambda)=0$.

\begin{corollary}
A radical square zero algebra is $\tau$-tilting infinite if and only if it admits a one-parameter family of bricks of the same length.
\end{corollary}

Before showing another property of min-$\tau$-infinite algebras, we remark that the center of a min-rep-inf algebra can be as large as the entire algebra. For instance, $k\langle x,y\rangle/\langle x^2, y^2, xy, yx \rangle$ is a commutative string algebra which is min-rep-inf.
This algebra is in fact obtained by gluing the source and sink in the Kronecker quiver. 
In contrast, our following result shows that the behaviour of min-$\tau$-infinite algebras with respect to their center is very different from their classical counterparts.

\begin{proposition}
Every minimal $\tau$-tilting infinite algebra is central.
\end{proposition}

\begin{proof}
As a consequence of \cite{EJR} and \cite{DIJ}, for any algebra $\Lambda$ and each ideal $J$ generated by some elements in $Z(\Lambda)\cap \rad(\Lambda)$, the sets $\texttt{i} \tau\rigid(\Lambda)$ and $\mathtt{i}\tau\rigid (\Lambda/J)$ are in bijection. In particular, $\Lambda$ is $\tau$-tilting finite if and only if $\Lambda/J$ is so. 

Now, assume $\Lambda=kQ/I$ is min-$\tau$-infinite. Obviously $|Q_0|=n
\geq 2$ (otherwise $\Lambda$ is local and therefore $\tau$-tilting finite).
Moreover, by the minimality assumption, we have $Z(\Lambda) \cap \rad(\Lambda)=0$.

Every element of $\Lambda=kQ/I$ can be expressed as $\rho:=\sum_{i=1}^{n} \lambda_i e_i+ r+I$, where $\lambda_i \in k$, and $r$ is a linear combination of paths of positive length in $Q$.
Suppose $\alpha$ is an arrow in $Q$ with $s(\alpha)=i$ and $e(\alpha)=j$. Then, multiplying $\rho$ by $\alpha$ from the left and from the right, and then reducing modulo $\rad^2 (\Lambda)$, we get $\lambda_i \alpha=\lambda_j \alpha$.
Since $Q$ is a connected quiver, we have all $\lambda_i$ are equal.
Consequently, we have $\rho= \lambda \cdot 1 + r + I$, for some $\lambda \in k$. Since $\lambda \cdot 1 + I \in Z(\Lambda)$, we get that $r + I \in Z(\Lambda)$. Moreover, from $Z(A) \cap \rad(\Lambda)= 0$ it follow that $r \in I$. Therefore, $\rho = \lambda \cdot 1$.
\end{proof}

As discussed in Section \ref{Subsection:Minimality conditions}, for an algebra $\Lambda$, both for the min-rep-inf and the min-$\tau$-infinite algebras the minimality conditions are defined with respect to certain sets of isormorphism classes of modules of finite length over $\Lambda$ and the corresponding sets on the quotient algebras $\Lambda/J$.
In our next proposition we give analogous characterizations of these minimality conditions in terms of certain modules of infinite lengths.
To do so, we recall an elegant characterization of representation infinite algebras by \cite{Au} Auslander. Recently, Sentieri \cite{Se} has shown a brick version of this theorem of Auslander. The following theorem collects both of the results.

\begin{theorem}[\cite{Au} and \cite{Se}]\label{Auslander-Sentieri}
Let $\Lambda$ be an algebra. Then,
\begin{enumerate}
    \item $\Lambda$ is representation finite if and only if any $M \in \Ind(\Lambda)$ is finitely generated.
    
    \item $\Lambda$ is brick finite if and only if any $M \in \Brick(\Lambda)$ is finitely generated.
    \end{enumerate}
\end{theorem}

Building upon the preceding theorem, we obtain new characterizations of the minimality conditions treated in this paper.

\begin{proposition}
Let $\Lambda$ be an algebra. Then,
\begin{enumerate}
    \item If $\Lambda$ is representation infinite, then it is minimal representation infinite if and only if every $M \in \Ind(\Lambda)\setminus \ind(\Lambda)$ is faithful.
    
    \item If $\Lambda$ is $\tau$-tilting infinite, then it is minimal $\tau$-tilting infinite if and only if every $N \in \Brick(\Lambda)\setminus \brick(\Lambda)$ is faithful.
    \end{enumerate}
\end{proposition}

\begin{proof}
To show $(1)$, suppose $\Lambda$ is min-rep-inf. If there exists $M \in \Ind(\Lambda)\setminus \ind(\Lambda)$ with $\ann_{\Lambda}(M)\neq 0$, then for the quotient algebra $\Lambda':=\Lambda/\langle \ann_{\Lambda}(M) \rangle$, we have $M \in \Ind(\Lambda')\setminus \ind(\Lambda')$. Since $\Lambda'$ is rep-finite, the desired contradiction follows from Theorem \ref{Auslander-Sentieri}.
    
    For the converse, let us assume that $\Lambda$ is not min-rep-infinite. Thus, for some nonzero ideal $J \in \Lambda$, the quotient algebra $\Lambda/J$ is rep-infinite and, again by Theorem \ref{Auslander-Sentieri}, there exists $M$ in $\Ind(\Lambda/J)\setminus \ind(\Lambda/J)$. But $M$ also belongs to $\Ind(\Lambda)\setminus \ind(\Lambda)$ and it is not faithful. So, we get the desired contradiction.

In the above proof, if we replace the indecomposable modules by bricks, we obtain a proof of $(2)$.
\end{proof}

\section{Minimal tilting infinite algebras}\label{Section:Minimal tilting infinite algebras}
Now we change our perspective and study the minimal $\tau$-tilting infinite algebras in the setting of classical tilting theory. Our main theorem in this section results in a good understanding of tilting modules and their mutation graph for minimal $\tau$-tilting infinite algebras. Before we show the main result of the section, we need some basic facts.

\begin{proposition}\cite[VIII.5.1]{ASS}\label{Prop:ASS on faithfulness}
A $\tau$-rigid $\Lambda$-module is partial tilting provided that it is faithful. In particular, $\tilt(\Lambda)$ consists of faithful modules in $\tau\ttilt(\Lambda)$.
\end{proposition}

For a basic $\Lambda$-module $M$, recall that $\Fac(M)$ denotes the subcategory of $\modu \Lambda$ which consists of all quotients of direct sums of $M$. In particular, if $M$ is $\tau$-rigid then $\Fac(M)$ is a functorially finite torsion class in $\modu \Lambda$. If $M$ belongs to $\taurigid (\Lambda) \setminus s\tau\ttilt(\Lambda)$, then it is a summand of a basic support $\tau$-tilting module $T$, which is given by the direct sum of all indecomposable $\Ext$-projective modules in $\Fac(M)$ (for details, see \cite{AIR}).

The next lemma is crucial in the proof of our main theorem of this subsection.

\begin{lemma} \label{LemmaTorsionAtBottom} Let $\Lambda$ be a minimal $\tau$-tilting infinite algebra. If $M$ is an unfaithful $\tau$-rigid module, then ${\rm Fac}(M)$ contains only finitely many bricks, up to isomorphism.
\end{lemma}

\begin{proof}
If $J:=\ann_{\Lambda}(M)$ then $JX=0$, for all $X$ in ${\rm Fac}(M)$. So, all bricks in ${\rm Fac}(M)$ are bricks over $\Lambda/J$. By Theorem \ref{Thm:tau-tilting finiteness equivalences} and the minimality of $\Lambda$, we are done.
\end{proof}

Now we are ready to prove the following theorem which establishes interesting connections between our results on min-$\tau$-infinite algebras and tilting theory. 

\begin{theorem}\label{Thm:almost all tau-rigid are faithful}
Let $\Lambda$ be minimal $\tau$-tilting infinite. Then, almost all $\tau$-rigid modules are faithful, thus partial tilting.
\end{theorem}
\begin{proof}
Suppose $\{T_i\}_{i\in \mathbb{N}}$ is a family of pairwise non-isomorphic unfaithful modules in $\taurigid(\Lambda)$, and let $\mathcal{T}_i = {\rm Fac}(T_i)$ be the corresponding functorially finite torsion classes.
From the beginning, we can assume every $T_i$ is a support $\tau$-tilting module.
By Lemma \ref{LemmaTorsionAtBottom}, each $\mathcal{T}_i$ has only finitely many bricks. Because every torsion class is fully determined by the bricks it contains, every interval $[0, \mathcal{T}_i]$ in the lattice $\tors(\Lambda)$ is finite. Since the Hasse quiver of the poset of functorially finite torsion classes of ${\rm mod}\Lambda$ is $n$-regular, we can find a strict ascending chain
$$0=\mathcal{U}_0 \subset \mathcal{U}_1 \subset \mathcal{U}_2 \subset \cdots$$
of functorially finite torsion classes with $\mathcal{U}_i \subseteq \mathcal{T}_i$, for each $i \in \mathbb{N}$. In particular, this is a chain of unfaithful torsion classes. If $J_i$ denotes the annihilator of $\mathcal{U}_i$, we get a corresponding descending chain
$$\cdots \subseteq J_2 \subseteq J_1 \subseteq J_0$$
of non-zero ideals. Since ideals of $\Lambda$ are finite dimensional, this chain has to stabilize to a non-zero ideal $I$, implying that there is a positive integer $r$ such that $I = J_{r+j}$ for all $j \ge 0$. For all $i \in \mathbb{N}$, observe that each $\mathcal{U}_i$, and therefore the algebra $\Lambda/J_i$, has at least $i$ non-isomorphic bricks. Thus, $\Lambda/I$ has infinitely many non-isomorphic bricks, and by Theorem \ref{Thm:tau-tilting finiteness equivalences}, $\Lambda/I$ is $\tau$-tilting infinite. This contradicts the minimality assumption on $\Lambda$. Thus, almost every $\tau$-rigid $\Lambda$-module is faithful.

The last assertion follows from the fact that, for any algebra, the family of faithful $\tau$-rigid modules coincides with that of partial tilting modules.

\end{proof}

\begin{remark}
For any min-$\tau$-infinite algebra $\Lambda$, Theorem \ref{Thm:almost all tau-rigid are faithful} implies that for almost all $\tau$-rigid $\Lambda$-modules we have $\pd_{\Lambda}(X)=1$, and therefore $\tau_{\Lambda}(X)=\rm{D}\Ext^1_{\Lambda}(X,\Lambda)$. This is interesting, particularly because there exist min-$\tau$-infinite algebras of infinite global dimension (see \cite{Mo2} for explicit examples).

In the previous theorem we note that the corresponding statement for bricks is not true. Namely, a min-$\tau$-infinite algebra may admit infinitely many isomorphism classes of unfaithful bricks (for instance, all regular bricks over the Kronecker algebra are unfaithful).
We further remark that there exist $\tau$-tilting infinite algebras $\Lambda$ such that almost all $\tau$-rigid $\Lambda$-module are faithful, but $\Lambda$ is not min-$\tau$-infinite (for example, consider the path algebra of the $3$-Kronecker quiver).
\end{remark}

Before we state an important consequence of Theorem \ref{Thm:almost all tau-rigid are faithful}, note that for an arbitrary algebra $\Lambda$, we may have $|\tau\ttilt(\Lambda) \setminus \tilt(\Lambda)|= \infty$. In particular, a $\tau$-tilting infinite algebra $\Lambda$ can admit only one tilting module.
Roughly speaking, our next result shows that for a min-$\tau$-infinite algebra $\Lambda$ the two sets $s\tau\ttilt(\Lambda)$ and $\tilt(\Lambda)$ are almost the same. This is far from obvious and suggests that min-$\tau$-infinite algebras can be useful in the study of some classical problems in tilting theory. In particular, it implies that over any min-$\tau$-infinite algebra $\Lambda$, the mutation graph of $\tilt(\Lambda)$ is infinite and $n$-regular at almost all vertices.

To state our result more precisely, we say $\Lambda$ is \emph{minimal tilting infinite} if it admits infinitely many tilting modules, up to isomorphism, but each proper quotient algebra $\Lambda/J$ has only finitely many tilting modules.

\begin{corollary}
If $\Lambda$ is a minimal $\tau$-tilting infinite algebra, then almost all support $\tau$-tilting modules are tilting.
Therefore, an algebra is minimal tilting infinite if and only if it is minimal $\tau$-tilting infinite.
\end{corollary}
\begin{proof}

This follows from the inclusions
$\tilt(\Lambda) \subseteq \tau\ttilt(\Lambda) \subseteq s\tau\ttilt(\Lambda)$,
and the fact that $\tau\ttilt(\Lambda)$ consists of sincere modules in $s\tau\ttilt(\Lambda)$.
By Theorem \ref{Thm:almost all tau-rigid are faithful}, almost all support $\tau$-tilting $\Lambda$-modules are faithful as $\Lambda$-modules, thus almost all $\tau$-tilting modules are tilting. 
It is evident that if $\Lambda$ is a min-$\tau$-inf algebra, then it must be minimal tilting infinite. If the converse fails, there exists a proper quotient $\Lambda/J$ which is min-$\tau$-infinite. Now, the desired contradiction immediately follows from the first part.
\end{proof}

\begin{remark}
For a given algebra $\Lambda$, it is \emph{a priori} a hard problem to describe the mutation graph of tilting modules, and even harder to decide whether all proper quotients of $\Lambda$ admit finite mutation graphs. This is partially because the notion of tilting finiteness is not preserved under taking quotients.
The previous corollary is significant in the sense that it yields an elegant classification of those algebras which are tilting infinite but all of their proper quotients are tilting finite.
Thanks to Theorem \ref{Thm:tau-tilting finiteness equivalences}, this classification could be stated in terms of bricks rather than tilting modules, which could have its own advantages, especially if one is interested in the study of geometric representation theory.
\end{remark}

From the preceding corollary and Theorem \ref{Thm: radical square zero}, the following result is immediate.

\begin{corollary}
Let $\Lambda=kQ/I$ with $\rad^2(\Lambda)=0$. If $\tilt(\Lambda)$ is an infinite set, then $Q$ contains a subquiver $Q'$ which is sink-source affine type.
\end{corollary}

\begin{remark}
It is well-known that the study of $\tau$-tilting theory closely relates to that of silting theory. In particular, for an algebra $\Lambda$, in \cite{AMV} the authors show that finite dimensional silting modules in $\Modu \Lambda$ are exactly the support $\tau$-tilting $\Lambda$-modules.
Thus, an algebra $\Lambda$ is $\tau$-tilting finite if and only if all torsion classes are generated by basic finite dimensional silting $\Lambda$-modules. This is the case if and only if all silting modules in $\Modu \Lambda$, up to some equivalence, are finite dimensional. This allows us to establish a connection between minimality conditions among $\tau$-tilting infinite and silting infinite algebras.
In particular, if $\Lambda$ is min-$\tau$-infinite, there exists a (definable) torsion class in $\Tors(\Lambda)$ such that the associated silting module, up to equivalence, is not finite dimensional.
It is immediate that if $\Lambda'$ is a proper quotient of $\Lambda$, all silting $\Lambda'$-modules, up to equivalence, are finite dimensional and hence $\Lambda'$ admits only finitely many basic silting modules. Therefore, as long as silting modules are considered up to equivalence, an algebra $\Lambda$ is min-$\tau$-infinite if and only if it is minimal silting infinite.
\end{remark}

\section{Reduction of the conjecture}\label{Section:Reduction of the conjecture}

In the current section we return to Conjecture \ref{main conjecture} and verify it for some families of minimal $\tau$-tilting infinite algebras. Consequently, we obtain Theorem \ref{Thm: reduction theorem in introduction}, as a reduction of Conjecture \ref{main conjecture} to a particular subfamily of min-$\tau$-infinite algebras.
Before we show the main result of this section, we need some preparation.

\begin{lemma}\label{Lemma:infinite family of ideals}
Let $\{I_i\}_{i \in \mathbb{N}}$ be an infinite subset of $\Ideal(\Lambda)$ such that the corresponding algebras $\Lambda/I_i$ are all isomorphic. If $\Lambda/I_1$ has a faithful brick of dimension $d$, then $\Lambda$ admits infinitely many non-isomorphic bricks of dimension $d$.
\end{lemma}

\begin{proof}
Let $B$ be a faithful brick of dimension $d$ in $\modu \Lambda/I_1$.
For each $i \ge 2$, let $B_i$ be the brick of $\Lambda/I_i$ induced from the isomorphism $\varphi_i: \Lambda/I_i \to \Lambda/I_1$. Then $B_i$ is faithful of dimension $d$. The $B_i$ are non-isomorphic as $\Lambda$-modules since they have pairwise distinct annihilators.
\end{proof}

Note that any algebra which satisfies the assumption of the previous lemma must be non-distributive. We postpone a more detailed study of $\tau$-tilting (in)finiteness of non-distributive algebras to our future work.

\medskip
Let $\Lambda$ be a minimal $\tau$-tilting infinite algebra such that it admits an infinite family $\{B_i\}_{i \in \mathbb{N}}$ of (non-isomorphic) unfaithful bricks, with $J_i:=\ann_{\Lambda}B_i$, for $i \in \mathbb{N}$. By the minimality assumption on $\Lambda$, these ideals must form an infinite family.
With no loss of generality, assume that they are pairwise distinct and consider the algebras $\Lambda/I_i$ for $i \ge 1$. Looking at these algebras up to isomorphisms, we have either a finite family, or an infinite family. In the former case, we apply Lemma \ref{Lemma:infinite family of ideals} to get an infinite family of non-isomorphic bricks of the same dimension. In case the $\Lambda/I_i$ form infinitely many non-isomorphic algebras, we have infinitely many $\tau$-tilting finite quotient algebras of the same dimension. 

Recall that if $\mathcal{T}$ and $\mathcal{T}'$ are functorially finite torsion classes in ${\rm mod}\Lambda$, such that $\mathcal{T} \subsetneq \mathcal{T}'$ is a covering relation in the poset of torsion classes, then there is a unique brick $B \in \mathcal{T}'$ which satisfies the following properties
\begin{enumerate}
    \item ${\rm Hom}_{\Lambda}(\mathcal{T},B)=0$.
    \item Every proper quotient of $B$ lies in $\mathcal{T}$.
    \item Every non-split short exact sequence $0 \to B \to X \to T \to 0$ with $T \in \mathcal{T}$ is such that $X \in \mathcal{T}$.
\end{enumerate}
Such a module is called a \emph{minimal extending brick} for $\mathcal{T}$. This brick is used to label the edge corresponding to the covering relation $\mathcal{T} \subsetneq \mathcal{T}'$ in the Hasse diagram of functorially finite torsion classes of $\modu \Lambda$. In \cite{BCZ} the authors extensively study the brick labelling of the poset of functorially finite torsion classes.

\medskip

In the poset of functorially finite torsion classes of $\Lambda$, a torsion class $\mathcal{T}$ is said to be at \emph{level} at most $t \in \mathbb{N}$  if the interval $[0, \mathcal{T}]$ contains a path of covering relations of length at most $t$.
For a functorially finite torsion class $\mathcal{T}$, in the next lemma we show that the length of all extending bricks for $\mathcal{T}$ are bounded by a function which depends only on the level of $\mathcal{T}$ and $d:=\dim_k \Lambda$. 
To give a more explicit proof for the next lemma, we use a bijection from \cite[Section 3]{AIR} which relates $s\tau \ttilt(\Lambda)$ and the isomorphism classes of basic $2$-term silting complexes in $K^b({\rm proj}\Lambda)$.
As further shown in \cite{AIR}, the aforementioned bijection gives rise to a poset isomorphism between the poset of support $\tau$-tilting modules and the poset of $2$-term silting complexes coming from \cite{AI}. Hence, the mutation of modules in $s\tau\ttilt(\Lambda)$ amounts to mutation of the corresponding $2$-term silting complexes, where the latter is in terms of certain triangles in the bounded derived category of $\modu \Lambda$. 
Namely, if $T < T'$ is a covering relation in the poset of $2$-term silting complexes in $K^b({\rm proj}\Lambda)$, by \cite[Theorem 2.31]{AI}, $T'$ is a right mutation of $T$ at one of its indecomposable summand.
More precisely, we have $T=U \oplus X$ and $T' = U \oplus Y$, with $U$ being a $2$-term presilting object and $X$ and $Y$ a pair of non-isomorphic indecomposable objects which appear in a triangle
$$Y \to U' \to X \to Y[1]$$
where $U' \to X$ is a minimal right ${\rm add}(U)$-approximation of $X$ (For full details, see \cite{AI} and \cite{AIR}.)

\begin{lemma}\label{lemma:bound on lenght of bricks}
Let $\mathcal{T}$ be a functorially finite torsion class in $\modu \Lambda$ which is at level at most $t$. Then, the dimension of any extending brick for $\mathcal{T}$ is bounded by a function that depends only on $t$ and on the dimension of $\Lambda$.
\end{lemma}

\begin{proof}
Let $T$ be a $2$-term silting complex in $K^b({\rm proj}\Lambda)$ represented by $P_1 \to P_0$. We assume that the differential is a radical morphism, hence $P_0$ and $P_1$ are uniquely determined up to isomorphism. Let $\dim T$ denote the total dimension of $T$, given by $\dim_k P_0+ \dim_k P_1$.

Using the same notation as in the paragraph preceding the lemma, we assume $T'$ is the right mutation of $T$ at $X$.
Let $\dim {\rm Hom}(U,X)=m$. Then we have a right ${\rm add}(U)$-approximation $U^m \to X$ of $X$. It is well known that there is a direct summand of $U^m$, isomorphic to $U'$, such that the corresponding restriction map is minimal. This yields $\dim U' \le m \, \dim U$. Therefore, 
\begin{eqnarray*}
\dim Y & \le & \dim U' + \dim X \\ & \le &  m \, \dim U + \dim X \\ & \le & (m+1) \dim T \\ & \le & (2(\dim T)^2 + 1) \dim T
\end{eqnarray*}
and hence $$\dim T' = \dim Y + \dim U \le (2(\dim T)^2 + 2)\dim T.$$

Note that for the minimal element of $s\tau\ttilt(\Lambda)$, the corresponding element in the poset of $2$-term silting complexes is $\Lambda \to 0$, which has total dimension $d:={\rm dim}_k\Lambda$.
Let $T$ be a $2$-term silting complex such that the interval $[0,T]$ in the Hasse quiver of poset of $2$-term silting complexes has a path of length at most $t$. Then there is a function $g(d,t)$ which depends only on $d$ and $t$ such that ${\rm dim}T \le g(d,t)$. Note that this function $g(d,t)$ is polynomial of degree $3^t$ in $d$.

Now, let $T$ be at level at most $t$ and $T'$ be a right mutation of $T$. By \cite[Theorem 3.2]{AIR}, the support $\tau$-tilting module $M'$ corresponding to $T'$ is the zeroth cohomology of $T'$. Thus, $\dim M' \leq \dim T' \le g(d,t+1)$.
By \cite[Theorem 2.7 and Theorem 3.2]{AIR}, we assume $\mathcal{T}$ denotes the functorially finite torsion class associated to $T$. Moreover, \cite[Corollary 2.34 and Corollary 3.9]{AIR} imply that $\mathcal{T}$ is at level at most $t$ in the poset of functorially finite torsion classes of $\modu \Lambda$, with the covering relation $\mathcal{T} < \Fac(M')$.
Since the minimal extending brick $B$ for the covering relation $\mathcal{T} < \Fac(M')$ is a quotient of $M'$, we obviously have $\dim B \leq g(d,t+1)$. This gives a function that bounds dimension of all minimal extending brick for $\mathcal{T}$.
\end{proof}

\begin{proposition}\label{Prop: Infinite family of pairwise distinct ideals}
Let $\Lambda$ be minimal $\tau$-tilting infinite and $\{I_i\}_{i\in A}$ an infinite family of pairwise distinct ideals in $\Lambda$. 
If each $\Lambda/I_i$ admits a faithful brick, then $\brick(\Lambda)$ has an infinite family of pairwise non-isomorphic bricks of the same dimension.
\end{proposition}

\begin{proof}
Without loss of generality, assume the family $\{I_i\}_{i\in A}$ is such that all ideals are of the same dimension, and this dimension is maximal with the property that each $\Lambda/I_i$ admits a faithful brick.
If infinitely many of these $\Lambda/I_i$ are isomorphic, apply Lemma \ref{Lemma:infinite family of ideals} and we are done. Thus, assume $\Lambda/I_i$ are pairwise non-isomorphic. If $V_i$ denotes a faithful brick over $\Lambda/I_i$, then $\{V_i\}_{i\in A}$ is a family of pairwise non-isomorphic modules in $\brick(\Lambda)$, particularly because $V_i$ have pairwise distinct annihilators in $\Lambda$. If there is a bound on the dimension of bricks in $\{V_i\}_{i\in A}$, we are done. Thus, for the sake of contradiction, we assume there is no bound on the dimension of $\{V_i\}_{i\in A}$.

For every functorially finite torsion class $\mathcal{T}$ at level $s$, by Lemma \ref{lemma:bound on lenght of bricks} there is a global bound on the dimension of all minimal extending bricks for $\mathcal{T}$.  Since all quotient algebras $\Lambda/I_i$ are $\tau$-tilting finite, there exists a positive integer $d_1$ such that infinitely many of the $\Lambda/I_i$ have a faithful brick in dimension $d_1$. For each $i$, by $\mathcal{H}_i$ we denote the Hasse quiver of the functorially finite torsion classes in $\modu \Lambda/I_i$. There is an arrow $\mathcal{T} \to \mathcal{T'}$ in $\mathcal{H}_i$ precisely when $\mathcal{T} < \mathcal{T}'$ is a covering relation.

Claim: For a fixed dimension $s$, we may assume that there is a function $h(s)$ such that every path in any given $\mathcal{H}_i$ has at most $h(s)$ minimal extending bricks of dimension at most $s$. 
To prove the claim, assume otherwise. Hence, there is a dimension $s' \le s$ such that we can find paths in the $\mathcal{H}_i$ having arbitrarily large number of minimal extending bricks of dimension $s'$. Note that in a given path, the minimal extending bricks are pairwise non-isomorphic. 
For the bricks of dimension $s'$ over all $\Lambda/I_i$, if their annihilators in $\Lambda$ form an infinite family of pairwise distinct ideals, then we are done. 
Therefore, we are in the case where there is an ideal $J$ such that $\Lambda/J$ has an infinite number of non-isomorphic bricks of dimension $s'$, and similarly we are done. This proves our claim.

By the claim, the number of bricks of dimension at most $d_1$ on any given path of each $\mathcal{H}_i$ is bounded by $h(d_1)$. Pick a positive integer $l_2 > h(d_1)$. We know that the $\mathcal{H}_i$  have infinitely many functorially finite torsion classes at level at most $l_2$. Therefore, by the pigeonhole principle and Lemma \ref{lemma:bound on lenght of bricks}, there is a positive integer $d_2 > d_1$ such that infinitely many of the $\Lambda/I_i$ have a brick of dimension $d_2$. In general, for $k \ge 2$, we pick $l_k > h(d_{k-1})$ and we get that there are infinitely many torsion classes of level at most $l_k$ in the $\mathcal{H}_i$. Hence, there is a positive integer $d_k > d_{k-1}$ such that infinitely many of the $\Lambda/I_i$ have a brick of dimension $d_k$.

 If for a given $r$, infinitely many of the bricks of dimension $d_r$ are faithful over the corresponding $\Lambda/I_i$, then we are done, as we obtain an infinite family of non-isomorphic bricks over $\Lambda$ of the same dimension $d_r$. 
Fix $r \ge 1$. Assume without loss of generality that each $\Lambda/I_i$ has an unfaithful brick $B_i$ of dimension $d_r$. Let $J_i:=\ann_{\Lambda}(B_i)$, then $I_i \subsetneq J_i$. If the ideals $J_i$ form an infinite family, then again, we get that the $B_i$ form an infinite family of non-isomorphic bricks of dimension $d_r$. We may thus assume that infinitely many of the $B_i$ have the same annihilator $J_{d_r}$. In particular, the algebra $\Lambda/J_{d_r}$ has a faithful brick of dimension $d_r$. We need only to consider the case where these $B_i$ are (almost all) isomorphic. 

Thus, we are left with the situation such that for each $r \ge 1$, we have an ideal $J_{d_r}$ with $\Lambda/J_{d_r}$ having a faithful brick of dimension $d_r$. We may assume that infinitely many of the ideals in the family $\{J_{d_r}\}$ are pairwise distinct, as otherwise a proper quotient of $\Lambda$ admits infinitely many unbounded bricks, and this contradicts the fact that each $\Lambda/J_{d_r}$ is $\tau$-tilting finite. Now, the family $\{J_{d_r}\}_{r \ge 1}$ yields an infinite family of ideals of the same dimension which are pairwise distinct and each $\Lambda/J_{d_r}$ admits a faithful brick. This gives the desired contradition with assumption on the dimension of ideals in the family $\{I_i\}_{i \in A}$. 
\end{proof}

From the previous proposition, we obtain the following theorem.

\begin{theorem}\label{Cor: Infinite many unfaithful bricks}
Let $\Lambda$ be a minimal $\tau$-tilting infinite algebra. If $\Lambda$ has infinitely many unfaithful bricks, then it admits an infinite family of bricks of the same length.
\end{theorem}

\medskip
Again, we note that the situation in Lemma \ref{Lemma:infinite family of ideals}, Proposition \ref{Prop: Infinite family of pairwise distinct ideals}, as well as the above theorem can occur only for non-distributive algebras. In particular, it is natural to ask whether it is true in general that a min-$\tau$-infinite algebra $\Lambda$ is non-distributive if and only if $\brick(\Lambda)$ contains an infinite family of unfaithful bricks. We remark that the analogous classification holds for min-rep-infinite algebras. Namely, a min-rep-inf algebra $\Lambda$ is non-distributive if and only if $\ind(\Lambda)$ contains an infinite family of unfaitfhul modules.
In our future work, we treat the non-distributive minimal $\tau$-tilting infinite algebras more closely.

\section{The First Brauer-Thrall conjecture for bricks}\label{Section: First BT conjecture}

In this short section, we use the geometric setting of representation varieties to give a new proof for the brick analogue of the First Brauer-Thrall conjecture, recently shown by Schroll and Treffinger \cite{ST}. To do so, we prove an interesting property of the minimal extending bricks in terms of the geometry of their orbits.

Recall that if the algebra $\Lambda$ is given by a bound quiver $(Q,I)$ over a field $k$, each module in $\modu \Lambda$ can be viewed as a point in a representation variety.
For a fixed dimension vector ${\bf{d}}=(d_x) \in \mathbb{Z}_{\geq 0}^{Q_0}$, we consider $\rep(Q,\bf{d})$ the affine space given by
$$\prod_{\alpha \in Q_1}{\rm Mat}_{\bf{d}_{e(\alpha)} \times \bf{d}_{s(\alpha)}}(k).$$
This variety parametrizes the representations of $kQ$ having dimension vector $\bf{d}$. More precisely, to the point $(V_{\alpha})_{\alpha \in Q_1}$ in $\rep(Q,\bf{d})$ is associated the representation $M$ with $M_x = k^{d_x}$ and $M(\alpha) = V_\alpha$.
Then, $\rep(Q,I,\bf{d})$ is the closed subset of $\rep(Q,\bf{d})$ whose points correspond to the representations annihilated by $I$. In this way, we view $\rep(Q,I,\bf{d})$ as an affine variety.
%Then, for $\Lambda=kQ/I$, points in the variety $\rep(Q,I,\bf{d})$ are in bijection with modules in %$\modu \Lambda$ whose dimension vector is $\bf{d}$.
Moreover, under the well-known action of the general linear group $\GL(\bf{d})$ on the variety $\rep(Q,I,\bf{d})$ via conjugation, the isomorphism classes of $\Lambda$-modules of dimension vector $\bf{d}$ are in bijection with the $\GL(\bf{d})$-orbits in $\rep(Q,I,\bf{d})$.
We sometimes denote $\rep(Q,I,\bf{d})$ by $\rep(\Lambda,\bf{d})$, where $\Lambda=kQ/I$.
We note that the variety $\rep(Q,I,\bf{d})$ is not necessarily irreducible, but has finitely many irreducible components and each such component is stable under the action of $\GL(\bf{d})$. 

\medskip

Before we state the following result, let us remark that for an arbitrary brick $M$ with $\underline{\dim}(M)= \bf{d}$, the $\GL(\bf{d})$-orbit of $M$ in $\rep(Q,I,\bf{d})$ is not necessarily open.
However, provided that $\Lambda$ is a brick finite algebras, every brick has a open orbit.
The next proposition gives us a nice geometric property of the minimal extending bricks defined in Section \ref{Section:Reduction of the conjecture}. In particular, this result is of significance over brick infinite algebras.

\begin{theorem}\label{MinExtBrickOpen}
Any minimal extending brick of a functorially finite torsion class has an open orbit.
\end{theorem}

\begin{proof}
Let ${\rm Fac}(N)$ be an arbitrary functorially finite torsion class, where $N$ is assumed to be a basic support $\tau$-tilting module corresponding to a $\tau$-rigid pair $(N,P)$ where $P$ is projective with $\Hom(P,N)=0$ and $N$ is $\tau$-tilting over its support. It follows from \cite[Proposition 4.13]{DI+} that any minimal extending brick $X$ labeling an arrow ${\rm Fac}(N) \to \mathcal{T}$ lies in $M^{\bot} \cap ^\perp\tau M \cap P^\perp$ where $M$ is a direct summand of $N$. By \cite[Theorem 3.8]{J}, the category $M^{\bot} \cap ^\perp\tau M \cap P^\perp$ is equivalent to the module category $\modu A$, for a finite dimensional local algebra $A$ (see also \cite[Theorem 4.12b]{DI+}).
%Let $(N,Q)$ be a $\tau$-rigid pair obtained from $(M,P)$ by removing exactly one indecomposable %summand, making it an almost complete $\tau$-rigid pair. 
%with $M=M_1\oplus M_2\oplus \cdots \oplus M_d$, and each $M_i$ is indecomposable, for all $1\leq i %\leq d$.
%Let $X$ be a minimal extending brick for ${\rm Fac}(M)$. Then, from \cite[Section 4.2]{DI+}, it %follows that there is a direct summand, say $M_1$, of $M$ such that for $N:=M_2\oplus \cdots \oplus %M_d$,
%we have that $X$ belongs to $N^{\bot} \cap ^\perp\tau N \cap P^\perp$ where $P$ is the maximal direct %summand of $\Lambda$ with ${\rm Hom}(P,M)=0$.
%By \cite[Theorem 3.8]{J}, the category $N^{\bot} \cap ^\perp\tau N \cap P^\perp$ is equivalent to the %module category $\modu A$, for a finite dimensional local algebra $A$.
Hence, $X$ is the unique brick in $M^{\bot} \cap ^\perp\tau M \cap P^\perp$.
Let $\underline{\dim}(X)=\bf{d}$ and consider the representation variety $\rep(\Lambda, \bf{d})$. 
Observe that the function $Y \mapsto {\rm dim}\Hom(M,Y)$ from $\rep(\Lambda, \bf{d})$ to $\mathbb{N}$ is upper-semicontinuous. Hence, the conditions $\Hom(M,-)=0$ defines an open set $U_1$ in $\rep(\Lambda, \bf{d})$. Similarly, the conditions, $\Hom(-,\tau M)=0$ and $\Hom(P,-)=0$ also define respective open sets $U_2$ and $U_3$ in $\rep(\Lambda, \bf{d})$. Therefore, $U:=U_1 \cap U_2 \cap U_3$ is a non-empty open set which contains $X$. Furthermore, let $\brick(\Lambda, \bf{d})$ denote the set of bricks in $\rep(\Lambda, \bf{d})$. Note that the orbit of each brick in $\rep(\Lambda, \bf{d})$ is of maximal dimension. Hence, it follows that $\brick(\Lambda, \bf{d})$ also forms an open set in $\rep(\Lambda, \bf{d})$. Consequently, we have a non-empty open set $U\cap \brick(\Lambda, \bf{d})$ in $\rep(\Lambda, \bf{d})$ which consists of bricks.
However, as mentioned above, $X$ is the unique brick in $M^{\bot} \cap ^\perp \tau M \cap P^\perp$, which implies that the orbit of $X$ must be open in $\rep(\Lambda, \bf{d})$.
\end{proof}

If a \emph{one-parameter family of bricks} means an irreducible curve of non-isomorphic bricks in the representation variety, the above proposition yields that none of the bricks in a one-parameter family of bricks can occur as a minimal extending brick. 

\medskip
 
Using the above proposition, we give an easy proof for the following theorem which could be viewed as a modern analogue of the first Brauer-Thrall conjecture. This result has been first shown in \cite{ST}.

\begin{theorem}
An algebra $\Lambda$ is $\tau$-tilting infinite if and only if there is no bound on the length of bricks in $\modu \Lambda$.
\end{theorem}

\begin{proof}
We only show the necessity, as the sufficiency follows directly from Theorem \ref{Thm:tau-tilting finiteness equivalences}.
%the fact that an algebra is $\tau$-tilting infinite if and only if it has infinitely many non-isomorphic bricks of finite length. Assume that $\Lambda$ is $\tau$-tilting infinite.
In the Hasse diagram of the functorially finite torsion classes of $\Lambda$, consider an infinite path starting at the minimal element. The corresponding minimal extending bricks form an infinite family of pairwise non-isomorphic bricks. We know that, for each dimension vector $\bf{d}$, the variety $\rep(\Lambda, \bf{d})$ has only finitely many irreducible components, and the closure of each open orbit provides an irreducible component. Hence, Proposition \ref{MinExtBrickOpen} implies that for each $\bf{d}$ there are only finitely many minimal extending bricks whose dimension vector is $\bf{d}$. 
Now the desired result is immediate.
\end{proof}

%As an immediate consequence of the previous theorem we observe that %the one-parameter family of bricks in the assertion of Conjecture %\ref{main conjecture} cannot be minimal extending bricks. %{\color{blue} Actually, no brick in a one-parameter family can be a %minimal extending brick.}{\color{purple} Can you explain why?
%A priori it seems we can have one minimal extending brick in an %irreducible component of $\rep(\Lambda, \underline{d})$ and an %infinite family of of the same length in the other irreducible %component!} {\color{blue} Right, but I think that a one-parameter %family is something that is entirely contained in a given %irreducible component.}

\medskip

\textbf{Acknowledgements.} The first-named author would like to thank Pierre-Guy Plamondon for his hospitality and several stimulating discussions during a $3$-month visit to Universit\'e du Paris Sud XI, supported by Mitacs Globalink Research Award. Moreover, he expresses gratitudes to the organizers of the trimester program on Representation Theory at Institut Henri Poincar\'e, held in 2020, where he developed part of this project.

\end{document}